\documentclass{article} \usepackage{amssymb} \newtheorem{theorem}{Theorem} \newtheorem{lemma}{Lemma}
 \newcommand{\La}{\Lambda}
\newcommand{\R}{{\mathbb R}}  \newcommand{\Z}{{\mathbb Z}} 
 \newcommand{\C}{{\mathbb C}}  

   \newcommand{\ef}{{\bf e}}

\begin{document} 
\begin{large}
\vspace{1cm}

\title{$L^2-$interpolation with error and size of spectra}

\author{Alexander Olevskii and Alexander Ulanovskii}

\date{} \maketitle

\begin{abstract}Given a compact set $S$ and a uniformly discrete sequence $\La$, we show that "approximate
interpolation"  of delta--functions on  $\La$ by a  bounded sequence of $L^2-$functions
with spectra in $S$ implies an estimate on measure of $S$ through the density of $\La$.
\end{abstract}

\section{ Introduction}

\medskip\noindent
Suppose $S$ is a bounded set on the real line $\R$. By $PW_S$ we shall denote the corresponding Paley--Wiener space:
$$
        PW_S:= \{f\in L^2(\R) ; \hat f=0 \mbox{ on } \R\setminus S\},
$$
where $$\hat f(t):=\frac{1}{\sqrt{2\pi}}\int_{\R}e^{itx}f(x)\,dx$$ denotes the Fourier transform.
It is well--known that each function $f\in PW_S$  can be extended to the
     complex plane as an entire function of finite exponential type.

Given a discrete set $\La=\{\lambda_j, j\in\Z\}\subset\R$, one says that $\La$ is a uniformly discrete  if $$
        \inf_{j\ne k} |\lambda_j-\lambda_k|   >0.
$$
This infimum  is  called the separation constant of $\La$.
 The following inequality  is well known (see
\cite{ya}, p. 82):
\begin{equation}\label{eque}
\Vert f\Vert_{L^2}\geq  C\Vert f|_\La\Vert_{l^2},
\ \mbox{for every} \ f\in PW_S.
\end{equation}
Here  $C>0$ is a constant which depends only on the separation constant of $\La$ and $S$,  $f|_\La$ denotes the restriction of $f$ on $\La$, and
$$
\Vert f\Vert^2_{L^2}:=\int_\R |f(x)|^2\,dx,\ \Vert f|_\La\Vert^2_{l^2}:=\sum_{j\in\Z}|f(\lambda_j)|^2.
$$
One can therefore regard the restriction $f|_\La$ as an element of $l^2(\Z)$,  the $j-$th coordinate of
$f|_\La$ being $f(\lambda_j)$.

\medskip\noindent
{\bf Definition}.  {\sl $\La$ is a  called a
{\sl set of interpolation} for $PW_S$, if for every data ${\bf c}\in l^2(\Z)$  there exists $f \in PW_S$ such that
\begin{equation}\label{int}
               f|_\La ={\bf c}.
\end{equation}}

A classical interpolation problem is to determine when $\La$ is a set of interpolation for $PW_S$.

 The upper uniform density of a uniformly discrete set $\La$ is defined as

 $$D^+(\Lambda):=\lim_{r\to\infty}\max_{a\in{\R}}\frac{\#(\Lambda\cap(a,a+r))}{r}.
$$
  A fundamental role of this quantity in the
interpolation problem, in the case when $S$ is
a single interval, was found by A. Beurling and J-P. Kahane.
 Kahane proved in \cite{k} that for $\La$ to be an
interpolation set for $PW_S$ it is necessary that
\begin{equation}\label{1}
         D^+(\La)\leq \frac{ 1}{2\pi}\mbox{ mes}\, S,
\end{equation}
and it is sufficient that
$$
         D^+(\La)< \frac{ 1}{2\pi}\mbox{ mes}\, S.
$$ Beurling (\cite{b1}) proved that  the last inequality  is necessary and sufficient
 for interpolation in the Bernstein
 space of all bounded on $\R$ functions with spectrum on the interval $S$.

The situation becomes much more delicate for
disconnected spectra, already when $S$ is a union of
two intervals. For the sufficiency part, not only the size but also the
arithmetics of $\La$ is important.
On the other hand, Landau \cite{l} extended the necessity part to the
general case:

\medskip\noindent{\bf Theorem A} {\sl
 Let $S$ be a bounded set. If a uniformly discrete set $\La$
is an interpolation set for $PW_S$ then  condition
(\ref{1}) is fulfilled.}

\section{Main result}
Denote by $\{\ef_j, j\in\Z\}$ the standard orthogonal basis in $l^2(\Z)$. When $S$ is compact, it is shown in \cite{ou1}  that Theorem A remains true under a weaker assumption that only $\ef_j, j\in\Z,$ admit interpolation by functions  from $PW_S$ whose norms are uniformly bounded.

Let us say that $\delta-$functions on $\La$ can be approximated  with error $d$ by functions from $PW_S$, if for every $j\in\Z$ there exists $f_j\in PW_S$ satisfying
\begin{equation}\label{vertj}\Vert f_j|_\La - \ef_j\Vert_{l^2} \leq d, \ j\in\Z.\end{equation}
The aim of this paper is to show that this `approximate'  interpolation of $\ef_j$ already gives an estimate on
the measure of $S$.
The result below extends both Theorem A (for compact $S$) and the mentioned result from \cite{ou1}.

\begin{theorem}\label{t5} Let $S$ be a compact set, and $\La$  a uniformly discrete set.
Suppose there exist functions $f_j\in PW_S$ satisfying (\ref{vertj}) for some $0<d<1$ and \begin{equation}\label{norms} \sup_{j\in\Z}\Vert f_j\Vert_{L^2}<\infty.\end{equation} Then
\begin{equation}\label{est}
                    D^+(\La)\leq \frac{1}{2 \pi(1- d^2)} \mbox{mes}\ S. \end{equation}
Bound (\ref{est}) is sharp for every $d$.
\end{theorem}

This result was announced in \cite{ou}.

Theorem 1 will be proved in sec. 4.
A variant of this result holds also when the norms of $f_j$ have a moderate growth, see sec. 5.

 
  \section{Lemmas}

\medskip\noindent 3.1. {\bf Concentration.}

   \medskip\noindent {\bf  Definition}: Given a number $c ,0<c<1,$ we say that a
linear subspace $X$ of $L^2(\R)$ is
 {\it $c$-concentrated} on a set $Q$ if $$
                     \int_Q |f(x)|^2\, dx \geq c \Vert f\Vert_{L^2}^2, \ f\in  X.
$$

\begin{lemma}\label{la} Given  sets $S,Q\subset\R$ of positive measure and a number $0<c<1$, let $X$ be a linear subspace of
$PW_{S}$ which is     $c$-concentrated on  $Q$.     Then
$$\mbox{dim}\, X \leq \frac{(\mbox{mes}\, Q)\,(\mbox{mes}\,
S)}{2\pi c}.$$
\end{lemma}

\medskip
    This lemma follows from  H.Landau's paper \cite{l} (see statements
(iii) and (iv) in Lemma 1, \cite{l}).

\medskip
\noindent
    3.2. {\bf A remark on Kolmogorov's width estimate.}

   \begin{lemma} Let $0<d<1$, and $\{{\bf u}_j\}, 1\leq j\leq n,$ be an orthonormal basis in an  $n$-dimensional complex
  Euclidean space $U$.
  Suppose that $\{{\bf v}_j\}, 1\leq j\leq n,$ is a family of vectors in $U$
  satisfying 
\begin{equation}\label{10}
                \Vert {\bf v}_j-{\bf u}_j\Vert \leq d,\ j=1,...,n.
\end{equation}
  Then for any $\alpha$, $
              1<\alpha < 1/d,$
 one can find a linear subspace $X$ in  $\C^n$ such that

         (i) dim$\, X >(1- \alpha ^2 d^2)n$-1,

         (ii) The estimate
$$
\Vert \sum_{j=1}^n c_j {\bf v}_j\Vert^2\geq (1-\frac{1}{\alpha})\sum_{j=1}^n|c_j|^2,
$$
holds for every vector $(c_1,...,c_n)\in X.$
\end{lemma}

\medskip
  \noindent
 The classical equality for Kolmogorov's width of "octahedron"
 (see \cite{ko}) implies that the dimension of the linear span of ${\bf v}_j$ is at least $(1-d^2)n$. This means that
  there exists a linear space $X\subset \C^n$, dim$\, X \geq (1-d^2)n$ such that
 the quadratic form
$$
                \Vert\sum_{j=1}^n c_j{\bf v}_j\Vert^2
$$
 is positive on the unite sphere of $X$.   The lemma above shows that by a
 small relative reduction  of the dimension, one can get an estimate of this
 form from below by a positive constant independent of $n.$

   We are indebted to E.Gluskin for the following simple proof of this
 lemma.

\medskip\noindent {\bf Proof}.
   Given an $n\times n$ matrix $T=(t_{k,l}), k,l=1,...,n$, denote by $s_1(T)\geq...\geq s_n(T)$
   the singular values of this matrix (=the positive square roots of the eigenvalues of $TT^*$).

   The following properties are well--known:

\medskip
        (a)  (Hilbert--Schmidt norm of $T$ via singular values)   $$\sum_{j=1}^n s_j^2(T) = \sum_{k,l=1}^n |t_{ k,l}|^2.$$

\medskip
        (b) (Minimax--principle for singular values)  $$s_k(T) = \max_{L_k}  \min_{x\in L_k,\Vert x\Vert=1}\Vert Tx\Vert,$$where the maximum is taken over all linear subspaces $L_k\subseteq \C^n$ of dimension $k.$

   \medskip     (c)  $s_{k+j} (T_1+T_2) \leq s_k(T_1)+s_j(T_2)$, for all $k+j\leq n.$

\medskip
   Denote by $T_1$  the matrix, whose columns are the coordinates of ${\bf v}_l$ in
   the basis ${\bf u}_k$, and   $T_2:=I-T_1$, where $I$ is the identity matrix. 
     Then  property (a) and (\ref{10}) imply: $$\sum_{j=1}^n s_j^2(T_2) < d^2 n,$$
      and hence:$$
                   s^2_j(T_2) \leq  d^2\frac{ n}{j},\ 1\leq j\leq n.
$$
   Now (c) gives:$$
           s_k(T_1) \geq 1-s_{n-k}(T_2)\geq 1- d\sqrt\frac{n}{n-k}  .
$$
  Taking the appropriate value of $k$, one can obtain  from (b)
  that there exists $X$ satisfying conclusions of the lemma.

\section{Proof of Theorem \ref{t5}}
1. Fix a small number $\delta>0$ and set $S(\delta):=S+[-\delta,\delta]$. Set
\begin{equation}\label{varp}g_j(x):=f_j(x)\varphi(x-\lambda_j),\
\varphi(x):=\left(\frac{\sin(\delta x/2)}{\delta x/2}\right)^2.
\end{equation}
 Clearly, $\varphi\in PW_{[-\delta,\delta]}$, so that $g_j\in PW_{S(\delta)}$. Also, since $\varphi(0)=0$ and $|\varphi(x)|\leq 1, x\in\R,$ it follows from (\ref{vertj})
that each $g_j|_\La$
approximates $\ef_j$ with an $l^2-$error $\leq d$:
\begin{equation}\label{gjl}
\Vert g_j|_\La-\ef_j\Vert_{l^2}\leq d, \ j\in\Z.
\end{equation}

\medskip\noindent 2. Given two numbers $a\in\R$ and $r>0$,  set $$  \La(a,r):=\La\bigcap (a-r,a+r),\ n(a,r):=\# \La(a,r).$$  For simplicity of presentation, in what follows we  assume that $\La(a,r)=\{\lambda_1,...,\lambda_{n(a,r)}\}$.

For every $g\in PW_{S(\delta)}$, we regard the restriction
$g|_{\La(a,r)}$ as vector in $\C^{n(a,r)}$.
It follows from (\ref{gjl}) that the vectors ${\bf v}_j:=g_j|_{\La(a,r)}$ satisfy (\ref{10}), where $\{{\bf u}_j, j=1,...,n(a,r)\}$ is the standard orthogonal basis
in $\C^{n(a,r)}$.

In the rest of this proof, we shall denote by $C$ different positive constants which do not depend on $r$ and $a$.

Fix a number $\alpha>1.$
 By Lemma 2,  there exists a subspace $X=X(a,r)\subset\C^{n(a,r)},$ dim$\,X\geq (1-\alpha^2 d^2)n(a,r)-1$, such that
$$
\Vert  \left(\sum_{j=1}^{n(a,r)} c_jg_j\right)|_{\La(a,r)}\Vert^2\geq (1-\frac{1}{\alpha})\sum_{j=1}^{n(a,r)} |c_j|^2, \ (c_1,...,c_{n(a,r)})\in X.
$$
By (\ref{eque}), this gives:
\begin{equation}\label{clr} \Vert\sum_{j=1}^{n(a,r)} c_jg_j\Vert_{L^2}^2\geq C\sum_{j=1}^{n(a,r)}|c_j|^2,\ (c_1,...,c_{n(a,r)})\in X.\end{equation}

\medskip\noindent 3. By  (\ref{norms}),  we have
$$
| f_j(x)|^2=\left|\frac{1}{\sqrt{2\pi}}\int_{S}e^{-itx}\hat f_j(t)\,dt\right|^2\leq
\frac{\mbox{mes}\,S}{2\pi}\Vert \hat f_j\Vert_{L^2}^2\leq C, \ j\in\Z.
$$Observe also that, since $\La$ is uniformly discrete, we have $n(a,r)\leq Cr$, for every $a\in\R$ and $r>1$.

\medskip\noindent 4.  Since $|x-\lambda_j|\geq \delta r$
whenever $\lambda_j\in(a-r,a+r)$ and $|x-a|\geq r+\delta r$, the inequalities in step 3 and (\ref{varp}) give
$$
\int_{|x-a|\geq r+\delta r}\left|\sum_{j=1}^{n(a,r)} c_jg_j(x)\right|^2\,dx=$$$$\int_{|x-a|\geq r+\delta r}\left|\sum_{j=1}^{n(a,r)}c_j
f_j(x)\left(\frac{\sin (\delta (x-\lambda_j)/2)}{\delta (x-\lambda_j)/2}\right)^2\right|^2\,dx \leq $$$$ Cr\left(\sum_{j=1}^{n(a,r)}|c_j|^2\right)
\int_{|x|>\delta r}\left(
\frac{2}{\delta x}\right)^4\,dx\leq \frac{C}{\delta^7 r^3}\left(\sum_{j=1}^{n(a,r)}|c_j|^2\right).
$$
This  and  (\ref{clr}) show that for every $\epsilon>0$ there exists $r(\epsilon) $ such that the linear space of functions
  $$g(x)=\sum_{j=1}^{n(a,r)}c_jg_j(x),\ (c_1,...,c_{n(a,r)})\in X,$$ is $(1-\epsilon)-$concentrated on $(a-r-\delta r,a+r+\delta r)$ for every $r\geq r(\epsilon)$, and every $a\in\R$.

\medskip\noindent
5. Lemma 1 now implies
   $$
\mbox{mes}\,S(\delta)\geq 2\pi(1-\epsilon)\frac{\mbox{dim}\,X}{\mbox{mes}\,(a-r-\delta r,a+r+\delta r)}$$$$\geq \frac{2\pi (1-\epsilon)}{1+\delta}  \frac{(1-\alpha^2 d^2)
\#\left(\La\bigcap (a-r,a+r)-1\right)}{2r}.
   $$
   Taking the limit as $r\to\infty$, where  $a=a(r)$ is such that the relative number of points of $\La$ in $(a-r,a+r)$ tends to $D^+(\La)$, we get
    $$
\mbox{mes}\,S(\delta)\geq \frac{2\pi (1-\epsilon)}{1+\delta}(1-\alpha^2 d^2)D^+(\La).
   $$
Since this is true for every $\epsilon>0,  \delta>0$ and  $\alpha>1$, we conclude that (\ref{est}) is true.

\medskip\noindent
Let us now check  that estimate (\ref{est}) in Theorem 1 is sharp.

\medskip\noindent{\bf Example}.  Pick up a number $a, 0<a<\pi$, and set $S:=[-a,a]$, $\La:=\Z$ and $$f_j(x):=\frac{\sin a(x-j)}{\pi (x-j)}\in PW_S, \ j\in\Z.$$ We have for every $j\in\Z$ that
$$
\Vert f_j|_\Z-{\bf e}_j \Vert^2_{l^2}=\Vert f_0|_\Z-{\bf e}_0
\Vert^2_{l^2}=\sum_{k\ne 0}\left(\frac{\sin a k}{\pi k}\right)^2
+\left(\frac{a}{\pi}-1\right)^2=
$$
$$
\frac{a}{\pi}-\frac{a^2}{\pi^2}
+\left(\frac{a}{\pi}-1\right)^2=1-\frac{a}{\pi}.
$$
Hence, the assumptions of Theorem 1 are fulfilled with
$d^2=1-a/\pi$. On the other hand, since $D^+(\Z)=1,$ we see that
mes$\, S= 2\pi(1-d^2)D^+(\Z)$, so that estimate (\ref{est}) is
sharp.

\section{Interpolation with moderate growth of norms}
When the norms of functions satisfying (\ref{vertj}) have a moderate growth
\begin{equation}\label{restr}\Vert f_j\Vert_{L^2}\leq Ce^{|j|^\gamma}, \ j\in\Z,\end{equation} where $C>0$ and  $0<\gamma<1$,    
the statement of Theorem 1 remains true, provided the density $D^+(\La)$ is replaced by the upper density $D^*(\La)$,
$$
D^{*}(\La):=\limsup_{a\to\infty}\frac{\#(\Lambda\cap(-a,a))}{2a}.
$$
Observe that $D^*(\La)\leq D^+(\La)$.

\begin{theorem}\label{t6} Let $S$ be a compact set, and $\La$  a uniformly discrete set.
Suppose there exist functions $f_j\in PW_S$ satisfying (\ref{vertj}) for some $0<d<1$ and (\ref{restr}). Then
 \begin{equation}\label{esti}
                   D^*(\La)\leq \frac{1}{2 \pi(1- d^2)} \mbox{mes}\ S. \end{equation}
 \end{theorem}

The upper density in this theorem cannot be replaced with the upper uniform density, see Theorem 2.5 in \cite{ou1}. The growth estimate (\ref{restr}) can be replaced with every `nonquasianalytic growth'. However, we do not know if it can be dropped.

\medskip\noindent  {\bf Proof of Theorem \ref{t6}}.
 The  proof  is  similar to the proof of Theorem~1.

\medskip\noindent
1. Fix  numbers $\delta>0$ and $\beta,$ $\gamma<\beta<1$. There exists a function $\psi\in PW_{(-\delta,\delta)}$ with the properties:
\begin{equation}\label{ee}
\psi(0)=1,\ |\psi(x)|\leq 1, \ |\psi(x)|\leq C e^{-|x|^\beta},\ x\in\R,
\end{equation}
where $C>0$ is some constant.
Such a function can be constructed as a product of $\sin(\delta_j x)/(\delta_jx)$ for certain sequence of $\delta_j$ (see Lemma~2.3 in \cite{ou1}).

Set$$h_j(x):=f_j(x)\psi(x-\lambda_j), \ j\in\Z.$$
Then each $h_j|_\La$ belongs to $PW_{S(\delta)}$ and
approximates $\ef_j$ with an $l^2-$error $\leq d$.

\medskip\noindent
2.   Set
$$
\Lambda_r:=\La\bigcap (-r,r),
$$
and denote by $C$ different positive constants independent on $r$.

 The argument in step 2  of the previous proof shows that there exists a linear space $X=X(r)$ of dimension $\geq (1-\alpha^2d^2)\# \La_r-1$
such that
$$
\Vert \sum_{j\in\La_r}c_jh_j(x)\Vert^2_{L^2}\geq C\sum_{j\in\La_r}|c_j|^2,
$$
for every vector $(c_j)\in X$.

\medskip\noindent
3.   Since $\La$ is uniformly discrete, we have $\#\La_r\leq Cr$ and $\max\{|j|, j\in\La_r\}\leq Cr, r>1$.  The latter estimate and (\ref{restr}) give:
$$
| f_j(x)|^2\leq\left(\frac{1}{\sqrt{2\pi}}\int_{S}\,|\hat f_j(t)|\,dt\right)^2\leq
\frac{\mbox{mes}\,S}{2\pi}\Vert \hat f_j\Vert_{L^2}^2\leq Ce^{Cr^\gamma}, \ j\in\Z.
$$

\medskip\noindent
4.  Using the estimates in step 3 and (\ref{ee}), we obtain:
$$
\int_{|x|\geq r+\delta r}\left|\sum_{j\in\Lambda_r} c_jh_j(x)\right|^2\,dx=$$$$\int_{|x|\geq r+\delta r}\left|\sum_{j\in\Lambda_r}c_j
f_j(x)\psi(x-\lambda_j)\right|^2\,dx \leq $$$$ \left(\sum_{j\in\Lambda_r}|c_j|^2\right)\left(Cr
e^{Cr^{\gamma}}\int_{|x|>\delta r}
e^{-2|x|^{\beta}}dx\right).
$$
Since $\beta>\gamma,$  the last factor  tends to zero as $r\to\infty.$ This and the estimate in step 2 show that  for every $\epsilon>0$ there exists $r(\epsilon) $ such that the linear space of functions
  $$\sum_{j\in\La_r}c_jh_j(x), \ (c_j)\in X,$$ is $(1-\epsilon)-$concentrated on $(-r-\delta r,r+\delta r)$, for all $r\geq r(\epsilon)$.

\medskip\noindent
5. Now, by Lemma 1, we obtain:
  $$
\mbox{mes}\,S(\delta)\geq  \frac{2\pi (1-\epsilon)}{1+\delta}  \frac{(1-\alpha^2 d^2)
\#\left(\La\bigcap (-r,r)-1\right)}{2r}.
   $$
By taking the upper limit as $r\to\infty,$ this gives    $$
\mbox{mes}\,S(\delta)\geq \frac{2\pi (1-\epsilon)}{1+\delta}(1-\alpha^2 d^2)D^*(\La).
   $$
Since this is true for every $\epsilon>0,  \delta>0$ and $\alpha>1$, we conclude that (\ref{esti}) is true.

\end{large}
\end{document}